\documentclass[12pt,a4paper]{article}
\usepackage{amsfonts}
\usepackage{amsmath}
\usepackage{mathrsfs}
\usepackage{graphicx}
\usepackage{amsmath,amsthm,amssymb,amscd}

\newtheorem{Thm}{Theorem}[section]
\newtheorem{Lem}[Thm]{Lemma}

\newtheorem{Def}[Thm] {Definition}

\newtheorem{Cor}[Thm]{Corollary}
\theoremstyle{remark}

\theoremstyle{claim}

\DeclareMathOperator{\Diff}{Diff}

\begin{document}

\begin{center}
{\Large \bf  Dominated Splitting and Pesin's Entropy Formula}\\
\smallskip
\end{center}

\bigskip
\begin{center}
Wenxiang Sun $^*$
\end{center}
\begin{center}
LMAM, School of Mathematical Sciences, Peking University, Beijing 100871, China\\
\end{center}
\begin{center}
E-mail: sunwx@math.pku.edu.cn
\end{center}
\smallskip
\begin{center}
Xueting Tian $^\dagger$
\end{center}
\begin{center} School of Mathematical Sciences, Peking University, Beijing 100871, China\\
\end{center}
\begin{center}
E-mail: txt@pku.edu.cn
\end{center}
\bigskip

\footnotetext {$^*$ Sun is supported by National Natural Science
Foundation ( \# 10671006, \# 10831003) and National Basic Research
Program of China(973 Program)(\# 2006CB805903) } \footnotetext
{$^\dagger$ Tian is the corresponding author.}
  \footnotetext {AMS Review: 37A05,
37A35, 37D25, 37D30}

\smallskip
\begin{abstract}
Let $M$ be a compact manifold and $f:\,M\rightarrow M$ be a $C^1$
diffeomorphism on $M$. If $\mu$ is an $f$-invariant probability
measure which is absolutely continuous relative to Lebesgue measure
and for $\mu$ $a.\,\,e.\,\,x\in M,$ there is a dominated splitting
$T_{orb(x)}M=E\oplus F$ on its orbit $orb(x)$, then we give an
estimation through Lyapunov characteristic exponents from below in
Pesin's entropy formula, i.e., the metric entropy $h_\mu(f)$
satisfies
$$h_{\mu}(f)\geq\int \chi(x)d\mu,$$ where
$\chi(x)=\sum_{i=1}^{dim\,F(x)}\lambda_i(x)$ and
$\lambda_1(x)\geq\lambda_2(x)\geq\cdots\geq\lambda_{dim\,M}(x)$ are
the Lyapunov exponents at $x$ with respect to $\mu.$ Consequently,
by using a dichotomy for generic volume-preserving diffeomorphism we
show  that Pesin's entropy formula holds for generic
volume-preserving diffeomorphisms, which  generalizes a result of
Tahzibi \cite{Ta} in dimension 2.
\end{abstract}

\section{Introduction}

To estimate    metric entropy through  Lyapunov exponents is an
important topic in differential ergodic theory. In 1977 Ruelle
\cite{Ru} got  from above an estimate  of metric entropy of an
invariant measure, and Pesin\cite{Pesin} in 1978 got  from below an
estimation  of metric entropy  of an  invariant  measure absolutely
continuous relative  to Lebesgue measure and thus got a so called
Pesin' entropy formula. Pesin's proof is based on the stable
manifold theorem. In 1980 Ma\~{n}\'{e} \cite{Ma} gave another
ingenious and very simple proof without using the theory of stable
manifolds. In 1985 Ledrappier and Young\cite{LY} generalized the
formula to all  SRB measures, not necessarily absolutely continuous
relative to Lebesgue measure. There are also more
generalizations\cite{Liu-endomorphism,Liu-random}.
\smallskip

Pesin's entropy formula by Pesin and by Ma\~{n}\'{e} and by others
assumes that not only the differentiability of the given dynamics is
of class $C^1$ but also that the first derivative satisfies an
$\alpha$-H$\ddot{o}$lder condition for some $\alpha>0$. It is
interesting to investigate Pesin's entropy formula under the weaker
$C^1$ differentiability hypothesis plus some additional condition,
for example, dominated splitting. The aim of this paper is to prove
that Pesin's entropy formula remains true for invariant probability
measure absolutely continuous relative to Lebesgue measure in the
$\textbf{C}^{\textbf{1}}$ diffeomorphisms with dominated splitting.
In the proof of \cite{Ma}, the combination of the graph transform
method (Lemma 3 there) and the distortion property deduced from the
H$\ddot{o}$lder condition of the derivative play important roles.
The domination assumption in  our $C^1$ diffeomorphism  helps us to
overcome much trouble. Our proof follows Ma\~{n}\'{e} without using
the theory of stable manifolds, as noted by  Katok that it seems
that Ma\~{n}\'{e}'s proof can also be extended to the more general
framework.
\smallskip

Tahzibi showed in \cite{Ta} that there is a residual subset
$\mathcal {R}$ in $C^{1}$ volume-preserving \textbf{surface}
diffeomorphisms such that every system in $\mathcal {R}$ satisfies
Pesin's entropy formula. As  an consequence  our main Theorem
\ref{PesFormula-Thm:1} and a result of Bochi and Viana\cite{BV}, we
generalize the result of Tahzibi into any dimensional case.

\bigskip

\section{Results}
Before stating our main results we need to introduce the concept of
dominated splitting. Denote the minimal norm of a linear map $A$ by
$m(A)=\|A^{-1}\|^{-1}$.

\begin{Def}\label{PesFormula-Def:2} Let $f: M\to M$ be a $C^1$ diffeomorphism on
a compact Remainnian manifold.

  (1). (Dominated splitting at one
point) Let $x\in M$ and $T_{orb(x)}M=E\oplus F$ be a $Df-$invariant
splitting on $orb(x)$. $T_{orb(x)}M=E\oplus F$ is called to be
$N(x)$-dominated at $x$, if there exists a constant $N(x)\in
\mathbb{Z}^+$ such that
$$\frac
{\|Df^{N(x)}|_{E(f^{j}(x))}\|}{m(Df^{N(x)}|_{F(f^{j}(x))})}\leq\frac12,\,\,\forall
\,j\in\mathbb{Z}.$$

(2). (Dominated splitting on an invariant  set) Let $\Delta$ be an
$f$-invariant set and $T_{\Delta}M=E\oplus F$ be a $Df-$invariant
splitting on $\Delta$. We call $T_{\Delta}M=E\oplus F$ to be a
$N$-dominated splitting, if there exists a constant $N\in
\mathbb{Z}^+$ such that
$$\frac
{\|Df^N|_{E(y)}\|}{m(Df^N|_{F(y)})}\leq\frac12,\,\,\forall \,y\in
\Delta.$$
\end{Def}

\bigskip

For a $Borel $ measurable map  $f: M\to M$ on a compact metric space
$M$ and an $ f-$invariant measure $\mu$, we denote by $h_\mu(f)$ the
metric entropy.

 Now we state our results as follows.
\begin{Thm}\label{PesFormula-Thm:1} Let $f: M\to M$ be a $C^1$ diffeomorphism on
a compact Remainnian manifold. Let $f$ preserve  an invariant
probability measure $\mu$ which is absolutely continuous relative to
Lebesgue measure. For $\mu$ a.e.  $x\in M,$ denote by
$$\lambda_1(x)\geq\lambda_2(x)\geq\cdots\geq\lambda_{dim\,M}(x)$$
the Lyapunov exponents at $x.$ Let $m(\cdot):M\rightarrow
\mathbb{N}$ be an $f$-invariant measurable function. If for
$\mu\,a.\,\,e.\,\,x\in M,$ there is a $m(x)$-dominated splitting:
$T_{orb(x)}M=E_{orb(x)}\oplus F_{orb(x)}$, then
$$h_{\mu}(f)\geq\int \chi(x)d\mu,$$ where
$\chi(x)=\sum_{i=1}^{dim\,F(x)}\lambda_i(x).$ \\In particular, if
for $\mu\,\, a.\,e.\,\,x\in M$, $E(x)$ and $F(x)$ coincide with the
sum of the Oseledec subbundles corresponding to negative Lyapunov
exponents and non-negative Lyapunov exponents respectively (or,
$E(x)$ corresponds to non-positive  Lyapunov exponents and $F(x)$
corresponds to positive Lyapunov exponents), then
$$h_{\mu}(f)=\int\chi(x)d\mu=\int\sum_{\lambda_i(x)\geq0}\lambda_i(x)d\mu.$$ In other words, Pesin's entropy
formula holds.
\end{Thm}

{\bf Remark.} Recall that the well known Ruelle's
inequality\cite{Ru}
$$h_{\mu}(f)\leq\int\sum_{\lambda_i(x)\geq0}\lambda_i(x)d\mu$$ is
valid for any invariant measure of $f$. Thus, if the inverse
inequality hold, the particular case of Theorem
\ref{PesFormula-Thm:1} is deduced immediately. So the left work we
need to prove is the inverse inequality.

\bigskip
Since Yang have proved in \cite{Yang} that for any diffeomorphism
$f$ far away from homoclinic tangency and any $f$-ergodic measure
$\mu$, the sum of the stable, center and unstable bundles in
Oseledec splitting is dominated on $\textrm{supp}(\mu),$ using
Theorem \ref{PesFormula-Thm:1} we have a direct corollary as
follows.
\begin{Cor}\label{PesFormula-Cor:2}
Let $f\in \Diff^1(M)$ far away from homoclinic tangency and let
$\mu$ be an $f$-ergodic probability measure which is absolutely
continuous relative to Lebesgue measure. Then $f$ satisfies Pesin's
entropy formula, i.e.,
$$h_{\mu}(f)=\sum_{\lambda_i>0}\lambda_i,$$ where
$\lambda_1\geq\lambda_2\geq\cdots\geq\lambda_{dim\,M}$ are the
Lyapunov exponents with respect to $\mu.$
\end{Cor}
\bigskip
 Let $m$ be the volume measure and $\Diff^1_m(M)$ denote the
space of volume-preserving diffeomorphisms. It is known that the
stable bundle and unstable bundle of Anosov diffeomorphism are
always dominated, and so are the bundles between the stable, center
and unstable directions in partially hyperbolic systems. Thus we
have a direct corollary as follows.

\begin{Cor}\label{PesFormula-Cor:3}Let $f\in \Diff^1_m(M).$ If $f$ is an Anosov
diffeomorphism (or, a partially hyperbolic diffeomorphism which
satisfies that for $m\,\,a.\,e.\,\,x,$ the Lyapunov exponents at $x$
in the central bundle are either all non-positive or all
non-negative), then Pesin's entropy formula holds.
\end{Cor}

\bigskip

In a Baire space, we say a set is residual if it contains a
countable intersection of dense open sets. We always call every
element in the residual set to be a generic point. It is known that
every $C^{1+\alpha}$ volume-preserving diffeomorphism satisfies
Pesin's entropy formula(see \cite{Ma,Pesin}) and the set of
$C^{1+\alpha}$ (or $C^2$) volume-preserving diffeomorphisms is dense
in $\Diff^1_m(M)$, so the set of volume-preserving diffeomorphisms
satisfying Pesin's entropy formula is dense in $\Diff^1_m(M)$.
Hence, it is natural to ask whether generic volume-preserving
diffeomorphisms satisfy Pesin's entropy formula. This problem is not
trivial because A. Tahzibi showed in \cite{Ta} that $C^{1+\alpha}$
volume-preserving diffeomorphisms are not generic in $\Diff^1_m(M)$.
Here we use Theorem \ref{PesFormula-Thm:1} to deduce this generic
property.

\begin{Thm}\label{PesFormula-Thm:2}
There exists a residual subset $\mathcal {R}\subseteq \Diff_m^1(M)$
such that for every $f\in\mathcal {R},$ the metric entropy
$h_{\mu}(f)$ satisfies Pesin's entropy formula, i.e.,
$$h_{\mu}(f)=\int\sum_{\lambda_i(x)\geq 0}\lambda_i(x)dm,$$ where
$\lambda_1(x)\geq\lambda_2(x)\geq\cdots\geq\lambda_{dim\,M}(x)$ are
the Lyapunov exponents of $x$ with respect to $m.$
\end{Thm}
{\bf Remark.} If $dim(M)=2,$ this result is firstly proved in
\cite{Ta}. \bigskip

\section{Proof of Theorem \ref{PesFormula-Thm:1}}

Our proof will be based on a general lower estimate for metric
entropy, which  makes it possible to avoid the use of partitions.
Let $g:M\rightarrow M$ be a map, $d$ be  a metric on $M$ and let
$\delta>0.$ If $x\in M$ and $n\geq0,$ define Bowen ball
$$B_n(g,\delta,x)=\{y\in M\,|\,d(g^j(x),g^j(y))\leq\delta,\,0\leq j\leq n\}.$$
In other words, $$B_n(g,\delta,x)=\bigcap_{j=0}^n g^{-j}
B_{\delta}(g^j(x)),$$ where $B_{\delta}(g^j(x))$ denotes the ball
centered at $x$ with  radius $\delta$. If  $g: M\to M$ is measurable
and $\mu$ is a measure on $M$(not necessarily $g-$invariant ),
define
$$h_{\mu}(g,\delta,x)=\limsup_{n\rightarrow+\infty}\frac1n[-\log\,\mu(B_n(g,\,\delta,\,x))].$$

\smallskip
\begin{Lem}\label{PesFormula-Lem:1} If $g$ is measurable, $\mu$ is a
$g$-invariant probability measure on $M$ and $\nu\gg\mu$ is another
measure on $M$ (not necessarily $g$-invariant), then
$$h_{\mu}(g)\geq\sup_{\delta>0}\int_M h_{\nu}(g,\delta,x)\,d\mu.$$
\end{Lem}

{\bf Proof} This lemma is a particular case of the Proposition in
\cite{Ma}, see P.96 or Lemma 13.4 in \cite{MaBook} for details.\hfill $\Box$ \bigskip

Before going into the proof of Pesin's formula we shall prove a
technical lemma. The reader familiar with the Hadamard graph
transform method for constructing invariant manifolds will recognize
this lemma one of the steps of that method. In the statement of the
lemma we shall use the following definitions from \cite{Ma,MaBook}.

\begin{Def}\label{PesFormula-Def:1}
Let $E$ be a normed space and $E=E_1\oplus E_2$ be a splitting.
Define $\gamma(E_1,E_2)$ as the supremum of the norms of the
projections $\pi_i:E\rightarrow E_i$ $i= 1, 2,$  associated with the
splitting. Moreover, we say that a subset $G\subset E$ is a $(E_1,
E_2)$-graph if there exists an open $U\subseteq E_2$ and a $C^1$ map
$\psi:U\rightarrow E_1$ satisfying $$G = \{x + \psi(x) |\,\,x \in
U\}.$$ The number $\sup\{\frac{\|\psi(x)-\psi(y)\|}{\|x
-y\|}|\,\,x\neq y \in U\}$ is called the dispersion of $G$.
\end{Def}
\bigskip

The following lemma about graph transform on dominated bundles is a
generalization to Lemma 3 in Ma\~{n}\'{e}\cite{Ma} about that on
hyperbolic bundles. Observe that the main  point   of the proof of
Lemma 3 there is the gap between two hyperbolic  bundles and can be
replaced by the gap of two dominated bundles, our proof of the
following lemma is a slight change of the proof of Lemma 3 in
Ma\~{n}\'{e}\cite{Ma}. We give a proof for completeness.

\bigskip

\begin{Lem}\label{PesFormula-Lem:2}
Given $\alpha>0,\,\,\beta>0$ and $c>0$, there exists $\tau>0$ with
the following property. If $E$ is a finite-dimensional normed space
and $E=E_1\oplus E_2$ a splitting with $\gamma(E_1,E_2)\leq\alpha,$
and $\mathcal{F}$ is a $C^1$ embedding of a ball $B_{\delta}(0)\subset E$ into
another Banach space $E'$ satisfying \\
$(i).\,\,\,\,\,\,\,\,\,\,\,$  $D_0\mathcal{F}$ is an isomorphism and
$\gamma((D_0\mathcal{F})E_1,(D_0\mathcal{F})E_2)\leq\alpha;$ \\
$(ii).\,\,\,\,\,\,\,\,\,$  $\|D_0\mathcal{F}-D_x\mathcal{F}\|\leq\tau$ for all $x\in B_{\delta}(0);$\\
$(iii).\,\,\,\,\,\,\,\,$ $\frac{\|D_0\mathcal{F}|_{E_1}\|}{m(D_0\mathcal{F}|_{E_2})}\leq\frac12;$\\
$(iv).\,\,\,\,\,\,\,\,\,$ ${m(D_0\mathcal{F}|_{E_2})}\geq\beta;$\\
then for every $(E_1,E_2)$-graph $G$ with dispersion $\leq c$
contained in the ball $B_{\delta}(0),$ its image $\mathcal{F}(G)$ is a
$((D_0\mathcal{F})E_1,(D_0\mathcal{F})E_2)$-graph with dispersion $\leq c. $
\end{Lem}

 {\bf Proof}
Identity $E$ with $E_1\times E_2$ and $E'$ with $(D_0\mathcal{F})E_1\times
(D_0\mathcal{F})E_2$. Write the map $F$ in the form
$$\mathcal{F}(x,y)=(Lx+p(x,y),\,\,Ty+q(x,y)),$$ where
$L=(D_0\mathcal{F})E_1,\,T=(D_0\mathcal{F})E_2.$ It follows that the partial derivatives
of $p$ and $q$ with respect to $x$ and $y$ have norm $\leq
\tau\alpha.$

Let $U\subset E_2$ be an open set and $\psi:U\rightarrow E_1$ a map
whose graph $\{(\psi(v),v)|v\in U\}$ is $G$. Then,
$$\mathcal{F}(G)=\{(L\psi(v)+p(\psi(v),v),\,Tv+q(\psi(v),v))|v\in U\}\}.$$ To
study this set define $\phi:U\rightarrow (D_0\mathcal{F})E_2$ by
$$\phi(v)=Tv+q(\psi(v),v)).$$ If $v,\,w\in U,$
$$\|\phi(v)-\phi(w)\|\geq\|T(v-w)\|-\|q(\psi(v),v)-q(\psi(w),w)\|.$$
Using the fact that the norm of the partial derivatives of $q$ are
$\leq\tau\alpha$ and hypothesis (iii) we obtain
$$\|\phi(v)-\phi(w)\|\geq m(T)\|v-w\|-\tau\alpha(\|\psi(v)-\psi(w)\|+\|v-w\|)$$
$$\geq(m(T)-\tau\alpha(1+c))\|v-w\|.\,\,\,\,\,\,\,\,\,\,$$
Hence, if $\tau$ is so small that $$m(T) - \tau\alpha ( 1 + c ) \geq
\beta-\tau\alpha(1+c)>0,$$ $\phi$ is a homeomorphism of $U$ onto
$\phi( U)$ whose inverse has Lipschitz constant $\leq (\beta -
\tau\alpha(1 + c))^{-1}$. In particular, $\phi( U)$ is open. Now
define $\hat{\psi}:\phi(U)\rightarrow (D_0\mathcal{F})E_1$ by
$$\hat{\psi}(v)=(L\psi\phi^{-1})(v)+p(\psi(\phi^{-1}(v)),\phi^{-1}(v)).$$ Clearly, $$\mathcal{F}(G) = \{(\hat{\psi}(x),
x)|x\in \phi(U).\} $$ To calculate the dispersion of $\mathcal{F}(G)$, write
$$\hat{\psi}=\tilde{\psi}\phi^{-1}$$ where $
\tilde{\psi}(w) = L\tilde{\psi}(w)+p(\tilde{\psi}(w), w).$ Then
$$\|\tilde{\psi}(w)-\tilde{\psi}(v)\|\leq \|L\|\|\psi(v)-\psi(w)\|+\tau\alpha(\|\psi(v)-\psi(w)\|+\|v-w\|)
$$$$\leq(c\|L\|+\tau\alpha(1+c))\|v-w\|.$$
 Then the dispersion of $\mathcal{F}(G)$ is less than or equal to $$c\frac{\|L\|+\tau\alpha(1+c)/c}{m(T)-\tau\alpha(1+c)}
 \,\,\,\,\leq\,\,\,\,
 c\frac{\frac12m(T)+\tau\alpha(1+c)/c}{m(T)-\tau\alpha(1+c)}$$
 $$=c\frac{\frac12+\tau\alpha(1+c)/cm(T)}{1-\tau\alpha(1+c)/m(T)}\,\leq\,
 c\frac{\frac12+\tau\alpha(1+c)/c\beta}{1-\tau\alpha(1+c)/\beta}.\,\,\,\,\,\,\,\,\,\,\,\,\,\,\,\,\,\,\,\,\,\,\,\,$$
 Taking $\tau$
small enough, the factor of $c$ is $<1$ and the lemma is
proved.\hfill $\Box$ \bigskip

\begin{Lem}\label{PesFormula-Lem:3}
Let $g\in \Diff^1(M)$ and $\Lambda$ be $g$-invariant subset of $M$.
If there is a $1$-dominated splitting on $\Lambda$:
$T_{\Lambda}M=E\oplus F$, then for any  $c>0$, there exists
$\delta>0$ such that for every $x\in \Lambda$ and any
$(E_x,F_x)$-graph $G$ with dispersion $\leq c$ contained in Bowen
ball $B_n(x,{\delta})\,\,(n\geq 0),$ its image $g^n(G)$ is a
$(D_xg^nE_x,D_xg^nF_x)$-graph with dispersion $\leq c. $
\end{Lem}

{\bf Proof} Let $\beta=\min_{x\in M}m(D_xg)$. Since dominated
splitting can be extended on the closure of $\Lambda$ and dominated
splitting is always continuous(see \cite{BLV}), we can take a finite
constant
$$\alpha=\sup_{x\in\Lambda}\gamma(E_x, F_x).$$  For given  $c>0$ and for the above $\alpha,\beta,$
take $\tau>0$ satisfying Lemma \ref{PesFormula-Lem:2}. Since $D_xg$
is uniformly continuous on $M,$ there is $\delta>0$ such that if
$d(x,y)<\delta,$ one has
$$\|D_xg-D_yg\|\leq\tau.$$  By applying
Lemma \ref{PesFormula-Lem:2}, we get the  following:\\
 \textbf{\,\,\,Fact\,\,\,}
 For any  $y\in\Lambda$ and every $(E_y,F_y)$-graph $H$ with dispersion $\leq c$
contained in the ball $B_{\delta}(y),$ its image $g(H)$ is a
$((D_yg)E_y,(D_yg)F_y)$-graph with dispersion $\leq c. $

We prove Lemma \ref{PesFormula-Lem:3} by induction. The conclusion
is trivial for $n=0$. Assume  it holds for some $n\geq0$, that is,
we assume that  if $G$ is a $(E_x,F_x)$-graph with dispersion $\leq
c$ contained in Bowen ball $B_{n}(g,{\delta},x)$ then  $g^n(G)$ is a
$(D_xg^{n}E_x,D_xg^{n}F_x)$-graph with dispersion $\leq c.$ Now let
$G$ is a $(E_x,F_x)$-graph with dispersion $\leq c$ contained in
Bowen ball $B_{n+1}(g,{\delta},x).$ Using
$B_{n+1}(g,{\delta},x)\subseteq B_{n}(g,{\delta},x),$ $G$ is also
contained in $ B_{n}(g,{\delta},x).$ So, by assumption $g^n(G)$ is a
$(D_xg^nE_x,D_xg^nF_x)$-graph with dispersion $\leq c$.
Take $y=g^n(x)\in\Lambda$ and let $H=g^n(G)$. Notice that $$(D_xg^nE_x,D_xg^nF_x)=(E_{g^nx},F_{g^nx})=(E_y,F_y)$$ and $$H=g^n(G)\subseteq g^n(B_{n}(g,{\delta},x))\subseteq
B_{\delta}(g^n(x))=B_{\delta}(y).$$
Thus $H$ is a
$(E_y,F_y)$-graph with dispersion $\leq c$ contained in $B_{\delta}(y)$.
Using the above \textbf{Fact}, we have $g(H)$
is a $((D_yg)E_y,(D_yg)F_y)$-graph with
dispersion $\leq c.$
Observe that $$g^{n+1}(G)=g(H)$$ and $$((D_xg^{n+1})E_x,(D_xg^{n+1})F_x)=((D_yg)E_y,(D_yg)F_y),$$ we get that $g^{n+1}(G)$
is a $((D_xg^{n+1})E_x,(D_xg^{n+1})F_x)$-graph with
dispersion $\leq c.$ \hfill $\Box$
\bigskip

\bigskip
Now we are ready to prove Pesin's formula.

{\bf Proof of Theorem \ref{PesFormula-Thm:1}} Put $$\Sigma_j=
\{x|\,\,\dim F(x)=j\}$$ and let $$S = \{j \geq 0 |\,\,\mu(\Sigma_j) >
0\}.$$ If $j\in S,$ let $\mu_j$ be the measure on $M$ given by
$$\mu_j(A) = \frac{\mu(A\cap\Sigma_j)}{\mu(\Sigma_j)}$$ for all
Borel subset $A$ of $M.$ Then $$\mu=\sum_{j\in S}\mu(\Sigma_j)\cdot
\mu_j$$ and thus by the affine property of metric entropy we have
$$h_{\mu}(f)=\sum_{j\in S}\mu(\Sigma_j)h_{\mu_j}(f).$$ Thus, all we
have to show is that
$$h_{\mu_j}(f)\geq\int\chi(x)d \mu_j.$$ This inequality obviously holds for $j=
0.$ Suppose $j> 0.$ Note that $\mu\ll Leb$ implies $\mu_j\ll Leb$
for all $j\in S.$ Hence, to simplify the notation we put
$$ \mu= \mu_j, \,\,\,\,\Sigma= \Sigma_j .$$

Fix any $\varepsilon > 0.$ Take $N_0$ so large that the set
$\Sigma_\varepsilon=\{x\in\Sigma|\,\,m(x)\leq N_0\}$ has
$\mu$-measure larger than $1-\varepsilon.$  Let $N=N_0!$ and
$g=f^N$, then the splitting $T_{\Sigma_\varepsilon}M=E\oplus F$
satisfies $1$-dominated with respect to $g$:
$$\frac{\|Dg|_{E(x)}\|}{m(Dg|_{F(x)})}\leq \prod_{j=0}^{\frac{N}{m(x)}-1}
\frac{\|Df^{m(x)}|_{E(f^{jm(x)}x)}\|}{m(Df^{m(x)}|_{F(f^{jm(x)}x)})}
\leq
(\frac12)^{\frac{N}{m(x)}}\leq\frac12,\,\,\,\forall\,x\in\Sigma_\varepsilon.$$
Note that $\Sigma_\varepsilon$ is $f$-invariant and thus
$g$-invariant. In what follows, in order to avoid a cumbersome and
conceptually unnecessary use of coordinate charts, we shall treat
$M$ as if it were a Euclidean space. The reader will observe that
all our arguments can be easily formalized by a completely
straightforward use of local coordinates.

Since dominated splitting can be extended on the closure of
$\Sigma_\varepsilon$ and dominated splitting is always
continuous(see \cite{BLV}), we can take and fix two constants $c>0$
and $a
> 0 $ so small that if $x\in\Sigma_\varepsilon\,, y \in M$ and $d(x, y)<a,$ then for every
linear subspace $E\subseteq T_yM$ which is a $(E(x), F(x))$-graph
with dispersion $<c$ we have $$\big{|}\log|detD_yg)|_E|-\log|det(D_xg)|_{F(x)}|\big{|}<\varepsilon.$$ Thus
\begin {equation}\label{PesFormula-Eq:1}
|detD_yg)|_E|\geq|det(D_xg)|_{F(x)}|\cdot e^{-\varepsilon}.
\end {equation} By Lemma \ref{PesFormula-Lem:3}, there exists
$\delta\in(0,\,a)$ such that for every $x\in \Sigma_\varepsilon$ and
any $(E_x,F_x)$-graph $G$ with dispersion $\leq c$ contained in the
ball $B_n(g,{\delta},x)\,\,(n\geq 0),$ its image $g^n(G)$ is a
$((D_xg^n)E_x,(D_xg^n)F_x)$-graph with dispersion $\leq c. $

 Let $\nu$ be the Lebesgue measure on
$M$. We give a claim as follows:\\
{\bf Claim.} For every $x\in \Sigma_\varepsilon$,
$$h_{\nu}(g,\delta,x)\geq N\chi(x)-\varepsilon.$$  By Lemma
\ref{PesFormula-Lem:1}, this property will imply that

$$\,\,\,\,\,\,\,\,\,\,\,\,\,\,\,\,\,\,\,\,\,\,\,\,\,\,\,\,\,\,\,\,\,\,\,\,\,\,\,\,\,\,\,
h_\mu(g)\geq\int_Mh_{\nu}(g,\delta,x)d\mu\,\,
\,\,\,\,\,\,\,\,\,\,\,\,\,\,\,\,\,\,\,\,\,\,\,\,\,\,\,\,\,\,\,\,\,\,\,\,\,\,\,\,\,\,
\,\,\,\,$$
$$\,\,\,\,\,\,\,\,\,\,\,\,\,\,\,\,\,\geq\,\,\,\,\,\int_{\Sigma_\varepsilon}h_{\nu}(g,\delta,x)d\mu\,\,
\,\,\,\,\,\,\,\,\,\,\,\,\,\,\,\,\,\,\,\,\,\,\,\,\,\,\,\,\,\,\,\,\,\,\,\,\,\,\,\,\,\,
$$$$\,\,\,\,\,\,\,\,
\,\,\,\,\,\,\,\,\,\,\geq
\,\,\,\,\,\int_{\Sigma_\varepsilon}(N\chi(x)-\varepsilon) d\mu\,\,
\,\,\,\,\,\,\,\,\,\,\,\,\,\,\,\,\,\,\,\,\,\,\,\,\,\,\,\,\,\,\,\,\,\,\,\,\,\,\,
$$
$$\,\,\,\,\,\,\,\,\,\,\,\,\,\,\,\,\,\,\,\,\,\,\,\,\,\,\,\,\,\,\,\,\,\,\,\,\,\,\,\,\,\,\, = \,\,\int_{M}N\,\,\chi(x) d\mu\,-\,\,
\int_{M\setminus \Sigma_\varepsilon}N\,\chi(x) d\mu\,-\,\varepsilon\cdot\mu(\Sigma_\varepsilon)
$$
$$\,\,\,\,\,\,\,\,\,\,\,\,\,\,\,\,\,\,\,\,\,\,\,\,\,\,\,\,\,\,\,\,\,\,\,\,
\,\,\,\,\,\,\,\,\,\,\,\,\,\,\,\,\,\,
\geq\,\,\,\, \int_{M}N\,\chi(x)d\mu-N\cdot C
\cdot dim(M)\cdot\mu(M\setminus \Sigma_\varepsilon)-\varepsilon$$
$$\,\,\,\,\,\,\,\,\,\,\,\,\,\,\,\,\,\,\,\,\,\,\,\,\,\,\,\,\,\,\,\,\,\,\,\,\,\geq \,\,\,\,
\int_{M}\,N\,\,\chi(x)\,d\mu\,-\,N \cdot C\cdot
dim(M)\cdot\varepsilon\,-\,\varepsilon$$ where $C=\max_{x\in
M}\log\|D_xf\|.$

 Hence, $$h_{\mu}(f)=\frac1Nh_{\mu}(g)\geq\int_{M}\chi(x)d\mu-C\cdot
dim(M)\cdot\varepsilon- \varepsilon. $$ Since $\varepsilon$ is arbitrary
this completes the proof of our theorem.

 It remains to prove the claim. Fix any $x\in\Sigma_\varepsilon.$ There exists $B > 0$ satisfying
$$\nu(B_n(g,\delta,x)) = B \int_{E(x)}\nu[(y+F(x))\cap
B_n(g,\delta,x)]d\nu(y)$$ for all $n \geq 0$, where $\nu$ also
denotes the Lebesgue measure in the subspaces $E(x)$ and $ y+F(x),\,y\in E(x).$  Thus
the claim is reduced to showing that
\begin {equation}\label{PesFormula-Eq:2}
\limsup_{n\rightarrow+\infty} \inf_{y\in E(x)}\frac1n[-\log
\nu(\Lambda_n(y))]\geq  N\chi(x)-\varepsilon,
\end {equation} where $$\Lambda_n(y)=(y+F(x))\cap
B_n(g,\delta,x).$$ If $\Lambda_n(y)$ is not empty, by Lemma
\ref{PesFormula-Lem:3} we have that $$g^n(\Lambda_n(y))\text{ is a }(E(g^n(x)),F(g^n(x)))\text{-graph with dispersion }\leq c.$$

Take $D>0$ such that $D>\textrm{vol} (G)$ (where $\textrm{vol}
(\cdot)$ denotes volume) for every $(E(w),F(w))$-graph $G$ with
dispersion $\leq c$ contained in
$B_\delta(w),\,w\in\Sigma_\varepsilon.$ Observe that $$g^n(\Lambda_n(y))\subseteq g^nB_n(g,\delta,x)\subseteq
B_\delta(g^n(x)),\,\,g^n(x)\in\Sigma_\varepsilon,$$ we have
$$D>\textrm{vol}(g^n(\Lambda_n(y)))=\int_{\Lambda_n(y)}|det(D_zg^n)|_{T_z\Lambda_n(y)}| d\nu(z).$$
Since $$g^j(\Lambda_n(y))\subseteq g^jB_n(g,\delta,x)\subseteq
B_\delta(g^j(x))\subseteq B_a(g^j(x)),\,\,j=0,1,2,\cdots,n,$$  we
have for any $z\in\Lambda_n(y)$,
$$d(g^j(z),g^j(x))<a,\,\,j=0,1,2,\cdots,n.$$ By inequality (\ref{PesFormula-Eq:1}), we
have
$$|det(D_zg^n)|_{T_z\Lambda_n(y)}|\,\,\,\,\,\,$$
$$=\,\,\prod_{j=0}^{n-1}|det(D_{g^j(z)}g)|_{T_{g^j(z)}g^j\Lambda_n(y)}|$$
$$\,\,\,\,\,\,\,\,\,\,\,\,\,\,\geq\,\,\prod_{j=0}^{n-1}\big{[}|det(D_{g^j(x)}g)|_{F(g^j(x))}|\cdot e^{-\varepsilon}\big{]}\,\,\,\,\,$$
$$=\,\,\,|det(D_{x}g^n)|_{F(x)}|\cdot e^{-n\varepsilon}.\,\,\,\,\,\,\,\,\,\,\,$$
Hence, $$\frac1n\log D\geq \frac1n\log\int_{
\Lambda_n(y)}|det(D_zg^n)|_{T_z\Lambda_n(y)}| d\nu(z)$$
$$\,\,\,\,\,\,\,\,\,\,\,\,\,\,\,\,\,\,
\,\,\,\,\,\,\,\,\,\,\,\,\,\, \geq\frac1n\log\int_{
\Lambda_n(y)}|det(D_{x}g^n)|_{F(x)}|\cdot e^{-n\varepsilon} d\nu(z)$$
$$\,\,\,\,\,\,\,\,\,\,\,\,\,\,\,\,\,\,
\,\,\,\,\,\,\,\,\,\,\,\,\,\, =\frac1n\log\big{[}\nu(
\Lambda_n(y))\cdot|det(D_{x}g^n)|_{F(x)}|\cdot e^{-n\varepsilon}\big{]}$$
$$\,\,\,\,\,\,\,\,\,\,\,\,\,\,\,\,\,\,\,\,\,\,\,\,
\,\,\,\,\,\,\,\,\,\,\,\,\,\, =\frac1n\log
\nu(\Lambda_(y))+\frac1n\log|det(D_{x}g^n)|_{F(x)}|-\varepsilon.$$
It follows that $$\lim_{n\rightarrow+\infty}-\frac1n\log
\nu(\Lambda_(y))\geq\lim_{n\rightarrow+\infty}\frac1n\log|det(D_{x}g^n)|_{F(x)}|-\varepsilon.$$
Combining this inequality and following equality from Oseledec
theorem\cite{Os}
$$\lim_{n\rightarrow+\infty}\frac1n\log|det(D_{x}g^n)|_{F(x)}|=N\chi(x),$$
we complete the proof of (\ref{PesFormula-Eq:2}). This completes the
proof of Theorem \ref{PesFormula-Thm:1}. \hfill $\Box$ \bigskip

\section{Proof of Theorem \ref{PesFormula-Thm:2}}
In this section we prove  Theorem \ref{PesFormula-Thm:2}. Before
that we need a result of Bochi and Viana\cite{BV}.
\begin{Thm}\label{PesFormula-Thm:3}(\cite{BV})
There is a residual subset $\mathcal {R}\subseteq \Diff_m^1(M)$ such
that for every $f\in\mathcal {R}$ and for $m\,a.\,\,e.\,\,x\in M$,
the Oseledec splitting of $f$ is either trivial(i.e., all Lyapunov
exponents are zero) or dominated at $x$.
\end{Thm}

\bigskip

{\bf Proof of Theorem \ref{PesFormula-Thm:2}} Let $\mathcal
{R}\subseteq \Diff_m^1(M)$ be the same as in Theorem
\ref{PesFormula-Thm:3}. Take and fix a diffeomorphism $f\in \mathcal
{R}.$ For $m\,\,a.\,e.\,\,x\in M,$ we can define
$$\chi(x)=\sum_{\lambda_i(x)\geq 0}\lambda_i(x).$$By Ruelle's
inequality\cite{Ru}, we have
$$h_{m}(f)\leq\int\chi(x)dm.$$ Thus we only need to prove that
$$h_{m}(f)\geq\int\chi(x)dm.$$

Let $$\Sigma_0=\{x\in M\,|\text{ the  Oselede  splitting  of  $f$
is trivial at }\,x\}$$ and $$\Sigma_1=\{x\in M\,|\text{ the  Oselede
splitting  of  $f$  is  dominated at }\,x\}.$$ Without loss of
generality, we assume that $m(\Sigma_0)>0$ and $m(\Sigma_1)>0$. Let
$m_j$ be the measure on $M$ given by
$$m_j(A) = \frac{m(A\cap\Sigma_j)}{m(\Sigma_j)}\,\,(j=0,1)$$ for all
Borel subset $A$ of $M.$ Then $m_0(\Sigma_0)=1,\,\,m_1(\Sigma_1)=1.$
More precisely, for $m_0\,a.\,\,e.\,\,x,$ the Oseledec splitting is
trivial at $x$ and for $m_1\,a.\,\,e.\,\,x,$ the Oseledec splitting
is dominated at $x$. Note that
$$m=m(\Sigma_0)\cdot m_0+m(\Sigma_1)\cdot m_1.$$ Thus by the affine property of metric entropy we have
$$h_m(f)=m(\Sigma_0)\cdot h_{m_0}(f)+m(\Sigma_1)\cdot h_{m_1}(f).$$
Based on these analysis we only need to prove that
$$h_{m_i}(f)\geq\int\chi(x)dm_i,\,i=0,1.$$

Since the metric entropy are always non-negative, obviously we have
$$h_{m_0}(f)\geq0=\int\chi(x)dm_0.$$
Note that $m_1$ are absolutely continuous relative to $m$. By
Theorem \ref{PesFormula-Thm:1}, we get
$$h_{m_1}(f)\geq\int\chi(x)dm_1.$$ This completes the proof of
Theorem \ref{PesFormula-Thm:2}. \hfill $\Box$

\bigskip

\bigskip

\section*{ References.}
\begin{enumerate}


\itemsep -2pt

\bibitem{BLV} Bonatti, Diaz, Viana, {\it Dynamics beyond uniform hyperbolicity: a global
geometric and probabilistic perspective}, Springer-Verlag Berlin
Heidelberg, 2005, 287-293.

\bibitem{BV} J. Bochi, M. Viana, {\it The Lyapunov exponents of generic volume preserving and symplectic
systems}, Ann. of Math., 161, 2005, 1423-1485.

\bibitem{LS} F. Ledrappier, J. Strelcyn, {\it A proof of the estimation from below in Pesin's entropy formula},
Ergod. Th. and Dynam. Sys., 2, 1982, 203-219.

\bibitem{LY} F. Ledrappier, L. S. Young, {\it The metric entropy of
diffeomorphisms}, Ann. of Math., 122,  1985, 509-539.

\bibitem{Liu-endomorphism} Peidong Liu, {\it Pesin's entropy formula for endomorphism,} Nagoya Math. J.
150 (1998) 197-209.

\bibitem{Liu-random} Peidong Liu, {\it Entropy formula of Pesin type for
non-invertible random dynamical systems,} Math. Z. 230 (1999) 201-39.

\bibitem{Ma} R. Ma\~{n}\'{e}, {\it A proof of Pesin's formula}, Ergod. Th. and Dynam. Sys., 1, 1981 95-102.

\bibitem{MaBook} R. Ma\~{n}\'{e}, {\it Ergodic theory and differentiable dynamics}, 1987, Springer-Verlag (Berlin, London).

\bibitem{Os} V. I. Oseledec, {\it Multiplicative ergodic theorem, Liapunov
characteristic numbers for dynamical systems}, Trans. Moscow Math.
Soc., 19, 1968, 197-221; translated from Russian.

\bibitem{Pesin} Y. Pesin, {\it Characteristic Lyapunov exponents and smooth ergodic
theory}, Russian Math. Surveys, 32, 1977, 55-114.

\bibitem{Ru} D. Ruelle,  {\it An inequality for the entropy of differentiable maps},
Bol. Sox. Bras. Mat, 9, 1978, 83-87.

\bibitem{Ta} A. Tahzibi, {\it $C^1$-generic Pesin's entropy formula
}, C. R. Acad. Sci. Paris, Ser. I 335, 2002, 1057-1062.

\bibitem{Yang} J. Yang, {\it $C^1$ dynamics far from tangencies}, preprint.
\end{enumerate}
\end{document}